\documentclass[11pt]{amsart}

\usepackage{amssymb}
\usepackage{amsmath}
\usepackage{amsthm}
\usepackage{amscd}
\usepackage{amsfonts}
\usepackage{ascmac}
\usepackage[usenames]{color}
\usepackage{graphics}

\theoremstyle{plain}
\newtheorem{thm}{Theorem}[section]
\newtheorem{prop}[thm]{Proposition}

\def\vol{\mathop{\mathrm{Vol}}\nolimits}

\newcommand{\bfC}{{\mathbf C}}
\newcommand{\bfP}{{\mathbf P}}
\newcommand{\bfR}{{\mathbf R}}

\newcommand{\barj}{{\overline j}}

\newcommand{\barpartial}{{\overline \partial}}

\newcommand{\mapright}[1]{\smash{\mathop{   \hbox to 0.7cm{\rightarrowfill}}
  \limits^{#1}}}

\def\bp{\overline{\partial}}

\def\grad{\mathrm{grad}}

\def\om{\omega}

\def\p{\partial}
\def\bp{\overline{\partial}}

\def\t{\theta}

\def\hh{\mathfrak h}

\newcommand{\Fut}{\mathrm{Fut}}
\newcommand{\Vol}{\mathrm{Vol}}

\renewcommand{\emph}[1]{{\color{red} \it #1}}

\definecolor{orange}{cmyk}{0, 0.7, 1, 0}
\definecolor{light-green}{cmyk}{0.5, 0, 0.5, 0}
\definecolor{light-blue}{cmyk}{0.5, 0, 0, 0}
\definecolor{light-yellow}{cmyk}{0,0,0.6,0}
\definecolor{dark-green}{cmyk}{0.7, 0, 0.7, 0.5}

\title{Volume minimization and obstructions to solving some problems in K\"ahler geometry}
\author{Akito Futaki and Hajime Ono}
\address{Graduate School of Mathematical Sciences, The University of Tokyo, 3-8-1 Komaba Meguro-ku Tokyo 153-8914, Japan}
\email{afutaki@ms.u-tokyo.ac.jp}
\address{Department of Mathematics, Saitama University, 255 Shimo-Okubo, Sakura-Ku,
Saitama 380-8570, Japan}
\email{hono@rimath.saitama-u.ac.jp}

\date{September 22, 2017}

\begin{document}
\begin{abstract} 
There is an obstruction to the existence of K\"ahler-Einstein metrics which is used to define
the GIT weight for K-stability, and it has been extended to various
geometric problems. This survey paper considers such extended obstructions to the existence problem of K\"ahler-Ricci solitons, Sasaki-Einstein metrics and (conformally) Einstein-Maxwell K\"ahler metrics. These three cases have a common feature that the obstructions are parametrized
by a space of vector fields.
We see, in these three cases, the obstructions are obtained as the derivative of suitable volume functionals. This tells us for which vector fields we should try to solve the existence problems.
\end{abstract}

\maketitle

\section{Introduction}
The existence of a K\"ahler-Einstein metric on a compact complex manifold $M$ has been known 
since 1970's in the case when $c_1(M) < 0$ by Aubin \cite{aubin76} and Yau \cite{yau78} where the K\"ahler class is the canonical class $K_M$, and in the case when $c_1(M) = 0$ by Yau \cite{yau78} where the K\"ahler class is arbitrary positive $(1,1)$-class.
In the remaining case when $c_1(M) > 0$, i.e. in the case when $M$ is a Fano manifold, 
the existence of a K\"ahler-Einstein metric is characterized by a condition called the K-stability by 
the recent works of 
Chen-Donaldson-Sun \cite{CDS3} and Tian \cite{Tian12}.
The K-stability is a condition in geometric invariant theory where the GIT weight, called 
the Donaldson-Futaki invariant \cite{donaldson02}, is defined extending an obstruction, now called 
the classical Futaki invariant,  obtained in \cite{futaki83.1}, \cite{futaki83.2}. The latter is defined for smooth compact
K\"ahler manifolds and is an obstruction to admit a
constant scalar curvature K\"ahler metrics (cscK metrics for short). Note that for a K\"ahler form
in the anti-canonical class on a Fano manifold, being a cscK metrics is equivalent to being a
K\"ahler-Einstein metric. On the other hand, 
the Donaldson-Futaki invariant is defined for possibly singular
central fibers of $\mathbf C^\ast$-equivariant degenerations, called the test configurations,
and a polarized K\"ahler manifold $(M,L)$ is said to be K-stable if the Donaldson-Futaki
invariant of the central fiber is non-negative for any test configurations and if the equality holds exactly when
the test configuration is product. 
Note that for the product configurations the Donaldson-Futaki invariant coincides with the classical
Futaki invariant.
The Fano case is the core of the conjecture known
as the Yau-Tian-Donaldson conjecture stating that a polarized K\"ahler manifold $(M,L)$ should admit
a cscK metric with its K\"ahler form in $c_1(L)$ if and only if 
$(M,L)$ is K-stable. In the K\"ahler-Einstein problem for the Fano case we take $L = K_M^{-1}$. 
The Yau-Tian-Donaldson conjecture for cscK problem with general polarizations is still remaining unsolved.
There are many variants of the Yau-Tian-Donaldson conjecture. For example, K-stability
characterizations for K\"ahler-Ricci solitons and Sasaki-Einstein metrics have been obtained 
respectively in \cite{DatarSzeke16} and \cite{CollinsSzeke15}. 
It is usually difficult to check whether a manifold is K-stable since there are infinitely 
many test configurations. However, in the cases with large symmetry groups checking K-stability
can be easier, see \cite{IltenSuss17}, \cite{Delcroix17}, \cite{Delcroix2016}. For alternate proofs for
the Yau-Tian-Donaldson conjecture for 
the Fano case, other important contributions, recent further developments and applications, the reader is referred to the two survey papers 
of Donaldson \cite{Donaldson2017a}, \cite{Donaldson2017b}.

The present survey paper focuses on extensions of the classical Futaki invariants for 
K\"ahler-Ricci solitons, Sasaki-Einstein metrics and Einstein-Maxwell K\"ahler metrics.
Existence problems for these three types of metrics have a common feature that they 
depend on the choice of a holomorphic
Killing vector field, and accordingly their obstructions have parameter space consisting
of holomorphic Killing vector fields in an appropriate Lie algebra.
The Ricci solitons are self-similar solutions of the Ricci flow and important object in the study
of singularity formations of the Ricci flow. On a compact K\"ahler manifold, 
a K\"ahler-Ricci soliton is a K\"ahler metric $g$ satisfying 
\begin{equation}\label{KRsoliton}
\mathrm{Ric}_g = g + L_{\mathrm{grad}f} g
\end{equation}
which is equivalent to
$$ \rho_g = \omega_g + i\partial\barpartial f$$
where $f$ is a Hamiltonian function for a holomorphic Killing vector field $X$, i.e. $X = J\mathrm{grad}f$, and $\rho_g$ and $\omega_g$ are respectively
the Ricci form and the K\"ahler form of $g$. Since $\rho_g/2\pi$ represents the first
Chern class, if a K\"ahler-Ricci soliton exists, the compact manifold $M$ is necessarily a Fano
manifold. Note also a Killing vector field on a compact K\"ahler manifold is necessarily holomorphic.
Given a Killing vector field $X$ we consider the toral group $T$ obtained by taking the closure of the flow generated by $X$, 
and ask
if there is a $T$-invariant K\"ahler-Ricci soliton $g$ satisfying (\ref{KRsoliton}) with
$X = J\mathrm{grad}f$. This problem is reduced to solving a Monge-Amp\`ere type equation.
However, Tian and Zhu \cite{TZ02} showed that there is an obstruction $\Fut_X$ to solving (\ref{KRsoliton}).
Thus if one choose an $X$ with non-vanishing $\Fut_X$ then one can never get a solution to the Monge-Amp\`ere equation. Tian and Zhu \cite{TZ02} showed that there is a twisted volume
functional $\Vol$ on the space of $X$ such that the derivative at $X$ of $\Vol$ is equal to $\Fut_X$:
\begin{equation}\label{derivative1}
d\Vol_X = \Fut_X.
\end{equation}
They further showed that the volume functional is proper and convex on the space of $X$.
Since holomorphic Killing vector fields on a compact K\"ahler manifold constitute a finite
dimensional vector space the volume functional has a unique minimum on the space of $X$. This gives the right choice to solve the equation (\ref{KRsoliton}). 

Sasaki-Einstein metrics caught considerable attention in mathematical physics through its 
 role in the AdS/CFT correspondence, and the volume minimization
is the key to find Sasaki-Einstein metrics. The Sasakian structure on an odd dimensional
manifold $S$ is by definition a Riemannian structure on $S$ such that its Riemannian cone
$C(S)$ has a K\"ahler structure. Fixing a complex structure on the cone $C(S)$, the deformation 
of the Sasakian structure on $S$ is given by the deformation of the cone structure on $C(S)$,
namely the deformation of the radial function $r$. The Reeb vector field is then given by
$Jr\partial/\partial r$. To each Reeb vector field one can assign a Sasakian structure on $S$.
Thus one can define the volume functional $\Vol$ on the space of Reeb vector fields.
The volume depends only on the Reeb vector field and is independent of the choice of the Sasakian structure with the given Reeb vector field. This fact is similar to the fact in K\"ahler geometry
that the volume depends only on the K\"ahler class and is independent of the choice of the
K\"ahler form in the given K\"ahler class.
The space of Reeb vector fields is the inside of 
the dual cone to the moment map image of the K\"ahler cone $C(S)$, and
the volume functional $\Vol$ is a homogeneous function on this space. 
Thus we may consider a slice which gives a
bounded domain sitting inside the dual cone. 
On the other hand the Sasaki-Einstein condition is equivalent to the 
K\"ahler cone $C(S)$ being Ricci-flat, 
and is also equivalent to the local transverse geometry of the Reeb flow being K\"ahler-Einstein
with positive scalar curvature.
One can then associate to each Reeb vector field $\xi$ an obstruction
$\Fut_\xi$ similarly to the Fano K\"ahler-Einstein problem \cite{FOW}, \cite{BGS}. 
Martelli-Sparks-Yau \cite{MSY2}
show for transversely Fano Sasakian manifolds
\begin{equation}\label{derivative2}
d\Vol_\xi = \Fut_\xi.
\end{equation}
In the case when $S$ is toric Sasakian, meaning when the cone $C(S)$ is toric K\"ahler, Martelli-Sparks-Yau further show that $\Vol$ is a proper convex function on the slice
in the dual cone consisting of the Reeb vector fields 
for which the volume functional $\Vol$ is defined. 
Thus there is a unique minimum $\xi$, and it is shown in \cite{FOW}, for any transversely Fano toric Sasakian manifold, there is a Sasaki-Einstein metric
with the choice of the unique minimum $\xi$ as the Sasakian structure.

Conformally K\"ahler Einstein-Maxwell metrics are relatively newer subject. The Einstein-Maxwell equation has been studied in general relativity in real dimension 4. In \cite{L1}, LeBrun
showed that, on a compact K\"ahler surface $(M,g)$, if there is a positive smooth function $f$ with
$J\mathrm{grad} f$ being a Killing vector field such that the Hermitian metric $\Tilde g = f^{-2}g$ 
has constant scalar curvature then $\Tilde g$ corresponds to a solution of the Einstein-Maxwell
equation. Thus, fixing a holomorphic Killing 
vector field $K$ and a K\"ahler class $\Omega$, to find a K\"ahler form $\omega_g \in \Omega$
such that $\Tilde g = f_K^{-2} g$ has constant scalar curvature is a problem in K\"ahler geometry, where $f_K$ is the Hamiltonian
function of $K$ with respect to $\omega$. In fact, if $K = 0$ then the problem is exactly the same
as the Yau-Tian-Donaldson conjecture.
Apostolov and Maschler \cite{AM} further set the problem into the Donaldson-Fujiki picture,
and formulated an extension $\Fut_K$ of the classical Futaki invariant parametrized by $K$.
In \cite{AM}, such $\Tilde g$ is called a conformally K\"ahler, Einstein-Maxwell metric. But we 
consider the problem of finding $(g,f_K)$ with $\omega_g$ in a fixed K\"ahler class, and therefore
it is more convenient to call such $g$ a (conformally) Einstein-Maxwell K\"ahler metric, or
even preferably omitting the word ``conformally''.
We then showed in \cite{FO17} that the derivative at $K$ of 
a suitably defined volume functional $\Vol$ on the space of $K$ satisfies
\begin{equation}\label{derivative3}
d\Vol_K = \Fut_K.
\end{equation}
However the volume functional is neither convex nor proper in general, and can have several
critical points.

In all these three cases, the critical points correspond to the cases when the classical Futaki
invariant vanishes. However, this may not be enough to have a solution, but the K-stability may be
the next issue. 

In section 2, 3 and 4 we give more details on 
K\"ahler-Ricci solitons, Einstein-Maxwell K\"ahler metrics and Sasaki-Einstein metrics respectively.

\vspace{0.3cm}

\noindent
Acknowledgment.\ 
The first auther would like to thank Yau Mathematical Sciences Center 
at Tsinghua University for its hospitality
where this survey paper was completed.

\section{K\"ahler-Ricci solitons}

In this section, we see how
a holomorphic Killing vector field
which admits a K\"ahler-Ricci soliton
is determined through the idea of volume minimization \cite{TZ02}.

Let $M$ be an $m$-dimensional Fano manifold.
A K\"ahler metric $g$ on $M$ with the K\"ahler form $\om_g\in
2\pi c_1(M)$ is called a {\it K\"ahler-Ricci soliton}
if there exists a holomorphic vector field $X$ on $M$ such that
\begin{equation}\label{eq:2.1}
\rho_g-\om_g=L_X\om_g
\end{equation}
holds, where $\rho_g$ denotes the Ricci form of $g$ and $L_X$ is the Lie
derivative along $X$.
In particular, if $X=0$, $g$ is a K\"ahler-Einstein metric.
Since $\rho_g$ and $\om_g$ represent $2\pi c_1(M)$,
there exists a real-valued smooth function $h_g$ such that
\begin{equation}\label{eq:2.2}
\rho_g-\om_g=i\p \bp h_g.
\end{equation}

On the other hand, for any holomorphic vector field $X$, the $(0,1)$-form
$\iota_X\om_g$ is $\bp$-closed. Therefore, by the Hodge theorem,
there exists a unique complex-valued smooth function $\t_X(g)$
such that
\begin{equation}\label{eq:2.3}
\iota_X\om_g=i\bp \t_X(g),\ \ \ 
\int_Me^{\t_X(g)}\om_g^m=\int_M\om_g^m.
\end{equation}

Hence we have
\begin{equation}\label{eq:2.4}
L_X\om_g=i\p \bp \t_X(g).
\end{equation}

By \eqref{eq:2.1}, \eqref{eq:2.2} and \eqref{eq:2.4},
a K\"ahler metric $g$ is a K\"ahler-Ricci soliton with respect to a
holomorphic vector field $X$ if and only if 
$h_g-\t_X(g)$ is constant.

It is difficult to determine $h_g-\t_X(g)$ explicitly.
However, Tian and Zhu \cite{TZ02} proved that
the integral of $v(h_g-\t_X(g))e^{\theta_X(g)}$ is independent of the
choice of $g$, where $v$ is a holomorphic vector  fields,
and it defines a holomorphic invariant.

\begin{thm}[\cite{TZ02}]\label{TZ-inv}
Let $\hh(M)$ be the Lie algebra which consists of all holomorphic vector fields on
$M$.  For a K\"ahler form $\om_g\in 2\pi c_1(M)$ and $X\in \hh(M)$, 
we define a linear function $\Fut_X$ on $\hh(M)$ as
\begin{equation}\label{eq:2.5}
\Fut_X(v)=\int_M v(h_g-\t_X(g))e^{\t_X(g)}\om_g^m,\ \ v\in \hh(M).
\end{equation}
Then $\Fut_X$ is independent of the choice of $\om_g\in 2\pi c_1(M)$.

If $M$ admits a K\"ahler-Ricci soliton with respect to $X\in \hh(M)$,
then $\Fut_X$ vanishes identically on $\hh(M)$.
\end{thm}

Note here that when $X=0$, this holomorphic invariant coincides with
the Futaki invariant, which is an obstruction to the existence of K\"ahler-Einstein
metrics in $c_1(M)$ \cite{futaki83.1}.

We next see that
the invariant $\Fut_X$ can be 
obtained as the first variation of some function on $\hh(M)$ \cite{TZ02}.
Such characterization of the holomorphic invariant plays a key role
in \S $3$ and \S $4$.

Let $X\in \hh(M)$. We renormalize the function $\t_X(g)$ defined by
\eqref{eq:2.3} to $\tilde{\t}_X(g)$ by adding a constant such that
\begin{equation}\label{eq:2.6}
\int_M\tilde{\t}_X(g)e^{h_g}\om_g^m=0.
\end{equation}

\begin{prop}[\cite{TZ02}]\label{1stvar}
Let a function $f$ on $\hh(M)$ be given by
\begin{equation}\label{eq:2.7}
f(Z)=\int_M e^{\tilde{\t}_Z(g)}\om_g^m.
\end{equation}
Then 
$f(Z)$ is independent of the choice of K\"ahler metrics with the K\"ahler class
$2\pi c_1(M)$.
Moreover the differential of $f$ at $X$ in the direction of $v\in \hh(M)$
is a constant multiple of $\Fut_X(v)$.
\end{prop}

By this proposition, if there exists a K\"ahler-Ricci soliton
with respect to a holomorphic vector field $X$, 
it is a critical point of $f$.

Let $\mathrm{Aut}^0(M)$ be the identity component of the
holomorphic automorphism group of $M$ and $K$ a maximal compact subgroup.
Then the Chevalley decomposition allows us to write
$\mathrm{Aut}^0(M)$ as a semi-direct product
\begin{equation}\label{eq:2.8}
\mathrm{Aut}^0(M)=\mathrm{Aut}_r(M)\ltimes R_u,
\end{equation}
where $\mathrm{Aut}_r(M)$ is a reductive algebraic subgroup of
$\mathrm{Aut}^0(M)$ which is the complexification of $K$, and
$R_u$ is the unipotent radical of $\mathrm{Aut}^0(M)$.
Let $\hh_r(M)$ and $\hh_u(M)$ be the Lie algebras of
$\mathrm{Aut}_r(M)$ and $R_u$ respectively. From the decomposition
\eqref{eq:2.8}, we obtain
\begin{equation}\label{eq:2.9}
\hh(M)=\hh_r(M)+\hh_u(M).
\end{equation}

\begin{prop}[\cite{TZ02}]\label{conv-proper}
Let $\vol$ be the restriction of $f$ to $\hh_r(M)$.
Then $\vol$ is a convex, proper real-valued function.
Hence there exists a unique minimum point $X_0\in \hh_r(M)$ of
$\vol$.
\end{prop}

By Proposition \ref{1stvar}, $\Fut_{X_0}$ vanishes identically on
$\hh_r(M)$. This minimum $X_0$ is the right choice to solve
the K\"ahler-Ricci soliton equation.
Note here that, to combine Proposition \ref{conv-proper} with
the result of Saito \cite{Saito},
$\Fut_{X_0}$ vanishes identically on $\hh(M)$.

For toric Fano manifold,
we can calculate $X_0$ as follows
\cite{Wang-Zhu}.
Let $M$ be an $m$-dimensional toric Fano manifold with the K\"ahler class
$c_1(M)$ and $\Delta_M\subset \bfR^m$ the 
corresponding moment polytope.
It is well-known that $\Delta_M$ is an
$m$-dimensional reflexive Delzant polytope.

Let $T$ be the maximal torus of $\mathrm{Aut}(M)$ and $\hh_0(M)$ 
its Lie algebra. 
$T$ is isomorphic to the $m$-dimensional algebraic torus $(\bfC^\times)^m$
and $\hh_0(M)$ is the maximal Abelian Lie subalgebra of $\hh(M)$.
If we take the affine logarithm coordinates
$(w_1,\dots,w_m)=(x_1+i\t_1,\dots,x_m+i\t_m)$
on $T\cong \bfR^m\times (S^1)^m$,
$\hh_0(M)$ is spanned by the basis
$\{\frac{\p}{\p w_1},\dots \frac{\p}{\p w_m}\}$.
Since $X_0\in \hh_0(M)$, $X_0$ can be expressed in the form
\begin{equation}\label{eq:2.10}
X_0=\sum _{i=1}^m c_i\frac{\p}{\p w_i}.
\end{equation}

\begin{prop}[\cite{Wang-Zhu}]\label{toric-HVF}
The constants $c_1,\dots,c_m$ in \eqref{eq:2.10} are given by the
following conditions:
\begin{equation}\label{eq:2.11}
\int_{\Delta_M}
y_i\exp \left\{\sum_{l=1}^m c_ly_l\right\}\,dy=0,\ \ i=1,\dots,m.
\end{equation}
\end{prop}

\section{Einstein-Maxwell K\"ahler geometry.}

In this section, we first introduce the notion
of conformally K\"ahler, Einstein-Maxwell (cKEM for short)
metrics defined by Apostolov-Maschler \cite{AM}
and give non K\"ahler examples
of cKEM metrics in any dimension.
We then define an obstruction to the existence of cKEM metrics
called cKEM-Futaki invariant and consider it from the view point
of volume minimization.
At the end of this section we give some results of computations 
on toric surfaces.

Let $(M,J)$ be a compact K\"ahler manifold. We call a Hermitian metric
$\Tilde{g}$ on $(M,J)$
a {\it conformally K\"ahler, Einstein-Maxwell metric}
 if it satisfies the following three conditions:

(a) There exists a positive smooth function $f$ on $M$  such that 
$g=f^2\Tilde{g}$  is K\"ahler.

(b) The Hamiltonian vector field $K=J\mathrm{grad}_gf$
 is Killing for both $g$ and $\Tilde{g}$.

(c)  $s_{\Tilde{g}}$ is constant.

As we mentioned in the Introduction, 
we call the K\"ahler metric $g$ in (a)
an {\it Einstein-Maxwell K\"ahler metric}.



By the definition above, cscK metrics are cKEM metrics.
However we consider them as trivial cKEM metrics.

The notion of cKEM metrics were introduced by
Apostolov-Maschler in \cite{AM}
as a generalization of strongly Hermitian solutions
of the Einstein-Maxwell equation.
We review some results by LeBrun on
strongly Hermitian solutions, see \cite{L1}, \cite{L2}.

Let $M$ be a compact manifold. 
A pair $(g,F)$ of a Riemannian metric $g$ and a real $2$-form $F$ is called
a solution of the
{\it Einstein-Maxwell equation} if it satisfies
$$
dF=0,\ d*_g F=0,\ [\mathrm{Ric}_g+F\circ F]_0=0,
$$
where $(F\circ F)_{jk} = 
F_j\,^\ell F_{\ell k}$ and $[\ ]_0$ denotes the trace free part.
This equation is the Euler-Lagrange equation of the following
functional which is studied in general relativity:
$$
(g,F)\mapsto \int_M \left(s_g+\lvert F\rvert^2_g \right) dv_g.
$$
LeBrun investigated Einstein-Maxwell equation
when $M$ is a complex surface in detail,
especially he introduced the notion of strongly Hermitian solutions:
Let $(g,F)$ be a solution of Einstein-Maxwell equation on a
complex surface $(M,J)$.
It is called a 
{\it strongly Hermitian solution} if
it satisfies
$$
\mathrm{Ric}_g(J\cdot,J\cdot)=\mathrm{Ric}_g(\cdot, \cdot),\ \ 
F(J\cdot,J\cdot)=F(\cdot,\cdot).
$$
LeBrun \cite{L1} pointed out that 
the metric component of a strongly Hermitian solution is a
cKEM metric.
Conversely, he also showed that for a cKEM metric $\tilde{g}$, one obtains
a strongly Hermitian solution
$$
(\tilde{g},\om_g+\frac12 f^{-2}[\rho_{\tilde{g}}]_0).
$$

We next give some examples
of cKEM metrics other than cscK metrics.
Typical known examples are conformally K\"ahler,
Einstein metrics by Page \cite{Page78} on the one point
blow up of $\bfC \bfP^2$, by Chen-LeBrun-Weber 
\cite{ChenLeBrunWeber} on the two point blow up
of $\bfC \bfP^2$, by Apostolov-Calderbank-Gauduchon \cite{ACG15}
on $4$-orbifolds and by B\'erard-Bergery \cite{BB82}
on $\bfC \bfP^1$-bundle over Fano K\"ahler-Einstein
manifolds.
Non-Einstein cKEM examples are constructed by LeBrun
\cite{L1}, \cite{L2} showing that there are ambitoric examples
on $\bfC \bfP^1\times \bfC \bfP^1$ and the one point
blow up of $\bfC \bfP^2$, and by Koca-T\o nnesen-Friedman \cite{KT}
on ruled surfaces of higher genus. The authors extended
LeBrun's construction on $\bfC \bfP^1\times \bfC \bfP^1$
to $\bfC \bfP^1\times M$ where $M$ is a compact
cscK manifold of arbitrary dimensions as follows \cite{FO17}.

Let $g_1$ be an $S^1$-invariant metric on $\bfC \bfP^1$ with
$\vol(\bfC \bfP^1,g_1)=2\pi$ and
$g_2$ a K\"ahler metric with $s_{g_2}=c$ on an $(m-1)$-dimensional compact complex manifold $M$.
The $S^1$-invariant metric $g_1$
can be written in the action-angle coordinates $(t,\theta)\in
(a,a+1)\times (0,2\pi]$ as
$$g_1=\frac{dt^2}{\Psi(t)}+\Psi(t)d\theta^2$$
for some smooth function $\Psi(t)$ which satisfies
the following boundary condition:
\begin{equation}\label{bc}
\Psi(a)=\Psi(a+1)=0,\ \ \Psi'(a)=-\Psi'(a+1)=2, \Psi>0\ \text{on }(a,a+1).
\end{equation}
Then we see that the constant scalar curvature equation 
$s_{\tilde{g}}=d$
for 
the metric $\tilde{g}=(g_1+g_2)/t^2$ on $\bfC \bfP^1\times M$
reduces to the following ODE:
\begin{equation}\label{eq:3.2}
t^2\Psi''-2(2m-1)t\Psi'+2m(2m-1)\Psi=ct^2-d.
\end{equation}

\begin{thm}[\cite{FO17}]
Let $c>8m-8$. Then there exist $a>0$ and $d>0$ such that
there exists a unique solution $\Psi$ of the ODE \eqref{eq:3.2}
which satisfies the condition \eqref{bc}.
As a result,  
for any  K\"ahler metric $g_2$ with $s_{g_2}=c$ on an $(m-1)$-dimensional compact complex manifold $M$, 
$$
\tilde{g}=\frac{1}{t^2}\left(\frac{dt^2}{\Psi(t)}+\Psi(t)d\theta^2+g_2\right)
$$
is an $S^1$-invariant cKEM metric on $\bfC \bfP^1\times M$.
\end{thm}

We next consider the existence problem of cKEM metrics.
Let $(M,J)$ be a compact complex manifold of $\dim_\bfC M=m$.
We fix a compact subgroup $G\subset \mathrm{Aut}(M,J)$,
a K\"ahler class $\Omega$, $K\in \mathfrak g$ and a 
sufficiently large $a\in \bfR$.
Denote by $\mathcal K^G_\Omega$ the space of $G$-invariant K\"ahler
metrics $g$ with $\om_g\in \Omega$. For $g\in \mathcal K^G_\Omega$,
there exists a unique function $f_{K,a,g}\in C^\infty(M)$ satisfying the following two
conditions:
\begin{equation}\label{eq:3.3}
\iota_K\omega_g=-df_{K,a,g},\ \int_M f_{K,a,g}\frac{\om_g^m}{m}=a.
\end{equation}
Note here that, for fixed $(K,a)$, $\min\{f_{K,a,g}(x)\,|\,x\in M \}$ is independent of 
$g\in \mathcal K^G_\Omega$, see \cite{AM}.
So if we choose $a$ sufficiently large, $f_{K,a,g}$ is positive for any 
$g\in \mathcal K^G_\Omega$.
Then we can ask the following existence problem;
does there exists a K\"ahler metric $g$ in $\mathcal K^G_\Omega$
such that $\tilde{g}_{K,a}=f^{-2}_{K,a,g}g$ is a cKEM metric?

When $K=0$, this is just the existence problem of
cscK metrics in $\mathcal K^G_\Omega$.
As a generalization of the Futaki invariant
\cite{futaki83.1}, \cite{futaki83.2},
Apostolov-Maschler \cite{AM} defined
the following integral invariant for non-zero $K$.

\begin{thm}[\cite{AM}]\label{cKEM-Futaki}
The linear function 
\begin{equation}\label{eq:3.4}
\Fut^G_{\Omega,K,a}:\mathfrak g\to \bfR,\ \ 
\Fut^G_{\Omega,K,a}(H):=
\int_M
\dfrac{s_{\tilde{g}_{K,a}}-c_{\Omega,K,a}}{f_{K,a,g}^{2m+1}}
f_{H,b,g}
\dfrac{\om_g^m}{m!},
\end{equation}
is independent of the choice of K\"ahler metric $g\in \mathcal K^G_\Omega$
and $b\in \bfR$. Here 
\begin{equation}\label{eq:3.5}
c_{\Omega,K,a}:=
\dfrac{\displaystyle{\int_M}s_{\tilde{g}_{K,a}}f_{K,a,g}^{-2m-1}\dfrac{\om_g^m}{m!}}
{\displaystyle{\int_M}f_{K,a,g}^{-2m-1}\dfrac{\om_g^m}{m!}}.
\end{equation}
is a constant which is independent of the choice of $g\in \mathcal K^G_\Omega$.
If there exists a K\"ahler metric $g\in \mathcal K^G_\Omega$
such that $\tilde{g}_{K,a}$ is a cKEM metric, then
$\Fut^G_{\Omega,K,a}$ vanishes identically.
\end{thm}

We call this linear function $\Fut^G_{\Omega,K,a}$ as the {\it cKEM-Futaki invariant}
for $(K,a)$.
We notice here that cKEM-Futaki invariant is parametrized by the pair
$(K,a)$. This situation bears resemblance to the holomorphic invariant
\eqref{eq:2.5} which is an obstruction to the existence of K\"ahler-Ricci solitons.
In fact, we now see that the cKEM-Futaki invariant can be
characterized as the first variation of the volume function.
To that end, we recall that constant scalar curvature Riemannian metrics
can be characterized as follows.
Let $M$ be a compact manifold with $n=\dim M\ge 3$ and 
$Riem(M)$ the set consists of all Riemannian metrics on $M$.
The scalar curvature $s_{g_0}$ of a Riemannian metric $g_0\in Riem(M)$
is constant if and only if $g_0$ is a critical point of the following
normalized Einstein-Hilbert functional on the conformal class of $g_0$:
\begin{equation}\label{eq:3.6}
EH(g):=\dfrac{\displaystyle{\int_M s_gdv_g}}{\displaystyle{\left(
\vol(M,g)\right)^{\frac{n-2}{n}}}}
\end{equation}

In our case, this functional gives the ``integral'' of the cKEM-Futaki invariant!

\begin{prop}[\cite{AM}]\label{inv-of-EH}
For a fixed $(K,a)$,
$EH(\tilde{g}_{K,a})$ is independent of the choice of $g\in \mathcal K^G_\Omega$.
\end{prop}

As a consequence, if there exists $g\in \mathcal K^G_\Omega$ such that
$\tilde{g}_{K,a}$ is a cKEM metric, then the pair $(K,a)$
is a critical point of the function
\begin{equation}\label{eq:3.7}
(K,a)\mapsto EH(K,a):=EH(\tilde{g}_{K,a}).
\end{equation}
The set of pairs
$$\mathcal P^G_\Omega
:=\{
(K,a)\in \mathfrak g\times \bfR\,|\,
f_{K,a,g}>0,\ g\in \mathcal K^G_\Omega
\}
$$
is a cone in the finite dimensional real vector space $\mathfrak g\times \bfR$.
Since the normalized Einstein-Hilbert functional is scale invariant,
the function $EH$ on $\mathcal P^G_\Omega$
reduces to the function on the quotient space $\mathcal P^G_\Omega/\bfR_+$.
If we choose representatives normalized as follows,
$EH$ can be represented as a power of the volume function.
We define a constant
$d_{\Omega.K,a}$ by
\begin{equation}\label{eq:3.8}
d_{\Omega,K,a}:=\dfrac{\displaystyle{\int_M s_{\tilde{g}_{K,a}}dv_{\tilde{g}_{K,a}}}}
{\vol(M,\tilde{g}_{K,a})}
=
\dfrac{\displaystyle{\int_M}s_{\tilde{g}_{K,a}}f_{K,a,g}^{-2m}\dfrac{\om_g^m}{m!}}
{\displaystyle{\int_M}f_{K,a,g}^{-2m}\dfrac{\om_g^m}{m!}}.
\end{equation}
By the argument in \cite{AM},
$d_{\Omega,K,a}$ is independent of the choice of $g\in \mathcal K^G_\Omega.$
Note here that, for general $(K,a)\in \mathcal P^G_\Omega$,
$c_{\Omega,K,a}\not=d_{\Omega,K,a}$. However if there exists a cKEM
metric $\tilde{g}_{K,a}$ then $c_{\Omega,K,a}=d_{\Omega,K,a}.$
Hence $c_{\Omega,K,a}-d_{\Omega,K,a}$ gives an obstruction to
the existence of cKEM metric $\tilde{g}_{K,a}$.
If we set
$$
\tilde{\mathcal P}^G_\Omega(\gamma):=
\{(K,a)\in \mathcal P^G_\Omega\,|\,
d_{\Omega,k,a}=\gamma\}
$$
for a constant $\gamma$, then
\begin{equation}\label{eq:3.9}
EH(K,a)=\gamma \vol (K,a)^{\frac{1}{m}}:=\gamma \vol (\tilde{g}_{K,a})^{\frac{1}{m}}
\end{equation}
on 
$\tilde{\mathcal P}^G_\Omega(\gamma)$.
By the first variation formula of the normalized Einstein-Hilbert functional
(cf. \cite{B}), we have
\begin{equation}\label{eq:3.10}
\frac{d}{dt}_{|t=0}EH(K+tH,a)=
\dfrac{2-2m}{\vol(K,a)^{\frac{m-1}{m}}}\int_M\left(\frac{s_{\tilde{g}_{K.a}}-d_{\Omega,K,a}}
{f_{K,a,g}
^{2m+1}}\right)f_{H,0,g}\frac{\om_g^m}{m!}
\end{equation}
and
\begin{equation}\label{eq:3.11}
\frac{d}{dt}_{|t=0}EH(K,a+tb)=\dfrac{2-2m}{\vol(K,a)^{\frac{m-1}{m}}}
(c_{\Omega,K,a}-d_{\Omega,K,a})\int_M
\frac{1}{f_{K,a,g}^{2m+1}}\frac{\om_g^m}{m!}.
\end{equation}

Therefore cKEM metrics have the following volume minimizing property.

\begin{thm}[\cite{FO17}]\label{cKEM-volmin}
Suppose that there exists a 
K\"ahler metric $g\in \mathcal K^G_\Omega$ such that
$\tilde{g}_{K,a}$ is a cKEM metric for 
$(K,a)\in \tilde{\mathcal {P}}^G_\Omega(\gamma)$.
Then $(K,a)$ is a critical point of the volume function
$\vol:\tilde{\mathcal P}^G_\Omega(\gamma)\to \bfR$.
Furter, $(K,a)$ is a critical point of $\vol$ if and only if
$\Fut^G_{\Omega,K,a}\equiv 0$.
\end{thm}

For example, 
let $(M,J,g)$ be an $m$-dimensional compact toric K\"ahler manifold.
We denote by $\Delta\subset \bfR^m$ the moment polytope.
Then we see that

\begin{equation}\label{eq:3.12}
EH(K,a)=
\frac{4\pi}{(m!)^{\frac{1}{m}}}
\frac
{\displaystyle{\int_{\partial \Delta}\frac{1}{f_{K,a}^{2m-2}}d\sigma}}
{\displaystyle{\left(\int_\Delta \frac{1}{f_{K,a}^{2m}}d\mu \right)^{\frac{m-1}{m}}}}
\end{equation}
for 
\begin{equation}\label{eq:3.13}
(K,a)\in \mathcal P^{T^m}_\Delta
\simeq
\{
f_{K,a}(\mu):=\sum_{i=1}^m K_i\mu_i+a\,|\,
f_{K,a}>0\text{ on }\Delta
\} 
\end{equation}
(cf. \cite{AM} or \cite{FO17}.)
Therefore, when $m=2$, we want to know
the critical points of
\begin{equation}\label{eq:3.14}
EH(a,b,c)^2=
8\pi^2
\frac
{\displaystyle{\left(\int_{\partial \Delta}\frac{1}{(a\mu_1+b\mu_2+c)^{2}}d\sigma\right)^2}}
{\displaystyle{\int_\Delta \frac{1}{(a\mu_1+b\mu_2+c)^{4}}d\mu }}
\end{equation}
with $a\mu_1+b\mu_2+c$ is positive on $\Delta$.
For $\bfC \bfP^2,\bfC \bfP^1\times \bfC \bfP^1$
and the one point blow up of $\bfC \bfP^2$,
we summarize results of computations.

\vspace{4mm}

$\bullet\ M=\bfC \bfP^2:$
In this case, up to scale and translations, $\Delta$ is the convex hull of the three points
$(0,0),(1,0)$ and $(0,1)$. The critical point of the function $EH$ on
$\mathcal P^{T^2}_\Delta/\bfR_+$
is only $[(0,0,1)]$. 

\vspace{3mm}

$\bullet\ M=\bfC \bfP^1\times \bfC \bfP^1:$
Let $\Delta_p$ be the convex hull of $(0,0),(p,0),(p,1)$ and $(0,1)$, where
$p\ge 1$.

When $1\le p\le 2$, $EH$ has the unique critical point
$[(0,0,1)]$.

On the other hand, when $p>2$, there exist three critical points
$$
[(0,0,1)],\ \left[\left(
\pm 1,0,\frac12\left(
\frac{p^{\frac32}}{\sqrt{p-2}}\mp p
\right)
\right)\right].
$$

We emphasize that this result shows that
the volume function is not convex unlike the case
of K\"ahler-Ricci solitons and of Sasaki-Einstein metrics, see \S $2$ and \S $4$.

\vspace{3mm}

$\bullet\ M=$ one point blow up of $\bfC \bfP^2:$

Let $\Delta _p$ be the convex hull of $(0,0),(p,0),(p,1-p)$ and
$(0,1)$, where $0<p<1$.

For $0<p<1$,
\begin{equation}\label{eq:3.15}
\left[
\left(
1,0,\frac{p(1-\sqrt{1-p})}{2\sqrt{1-p}+p-2}
\right)
\right]
\end{equation}
is a critical point of $EH$.

When $\frac89 <p<1$ there are the following two more critical points
\begin{equation}\label{eq:3.16}
\left[
\left(
-1,0,
\frac{p(3p\pm \sqrt{9p^2-8p})}{2(p\pm \sqrt{9p^2-8p})}
\right)
\right].
\end{equation}

Let $0<\alpha <\beta<1$ be the real roots of 
$$
F(p):=p^4-4p^3+16p^2-16p+4=0.
$$

When $0<p<\alpha$, there are the following two critical points

\begin{equation}\label{eq:3.17}
\left[
\left(
p^2-4p+2\pm \sqrt{F(p)},\pm 2\sqrt{F(p)},
p^2+2p-2\mp \sqrt{F(p)}
\right)
\right].
\end{equation}

An extension of Lichnerowicz-Matsushima theorem asserting the reductiveness of the
automorphism group on a cKEM manifold is obtained in \cite{FO_reductive17} and \cite{Lahdili17}.

\section{Sasakian Geometry.}

A Sasakian structure is often referred to as an odd dimensional analogue of the K\"ahler structure.
It roughly consists of a contact structure, a Riemannian structure compatible with the contact strucure 
and an almost complex structure on the
contact bundle.
There are many equivalent definitions, but the following definition is the most simple and rigorous one.
In Riemannian point view, a Sasakian manifold 
is a Riemannian manifold $(S, g)$ 
whose cone manifold 
$(C(S), \overline g)$ with $C(S) \cong  S\times \bfR^+$ and $\overline g = dr^2  + r^2g$ 
is K\"ahler, where $r$ is the standard coordinate on $\bfR^+$. 
In this paper we always assume $S$ is closed and connected. 
From the definition, $S$ is odd-dimensional and we put $\dim S = 2m + 1$, and thus
$\dim_\bfC {C(S)} = m+1$.  $S$ is identified with the submanifold $\{r=1\} \subset C(S)$. 
The K\"ahler form on $C(S)$ is given by $i\partial\barpartial r^2$. From this we see that,
fixing the holomorphic structure on $C(S)$, 
the Sasakian structure is determined by the radial function $r$ since the Riemannian structure
is induced from the K\"ahler structure of $C(S)$. 
We consider the deformations of the Sasakian structure on $S$
fixing the complex structure $J$ on $C(S)$.

On the other hand, $S$ also inherits a contact structure 
with the contact form
$$\eta = (i(\barpartial - \p) \log r)|_{r=1}.$$
It is well known \cite{BGbook} that the Sasakian structure is determined
by the transverse K\"ahler structure of the flow generated by the Reeb vector field $\xi$ of $\eta$.
The Reeb vector field $\xi$ is obtained by restricting the vector field $\tilde \xi := J(r\frac{\p}{\partial r})$ on $C(S)$ to
$S = \{r=1\} \subset C(S)$. 
This is a standard fact known as the ``K\"ahler sandwich'': The Sasakian structure is
equivalently given by the K\"ahler strcuture on the cone or given by the transverse
K\"ahler structure on the local orbit spaces of the Reeb flow, see \cite{BGbook} for the
detail.
From this we see that the Sasakian structure can be deformed
by the deformation of the choice of Reeb vector field in the Lie algebra $\mathrm{Lie}(T_\xi)$ of
the torus $T_\xi$ obtained by taking the closure of the flow generated by $\xi$ since the 
deformed Reeb flow still has transverse K\"ahler structure.
Then by choosing a rational point in $\mathrm{Lie}(T_\xi)$ we obtain a Reeb vector field obtained
as an $S^1$-action on an ample line bundle over an orbifold. Thus $C(S)$ has an affine 
algebraic variety $\mathcal A$ with only one singular point at the apex as an underlying space.

Let $G$ be the group of biholomorphisms of $\mathcal A = C(S)$ preserving the cone structure, that is,
$\mathrm{Lie}(G)$ consists of the real parts of holomorphic vector fields on $\mathcal A$ commuting with $r\frac{\p}{\partial r}$.
Let $T$ be the maximal torus of $G$ containing $T_\xi$. Note here that it is a standard fact that $r\frac{\p}{\partial r}$
preserves $J$ and that $\tilde \xi - iJ\tilde \xi$ is a holomorphic vector field. The deformation space of $T$-invariant Sasakian structures  containing the Sasakian structure of $S$,
or equivalently $T$-invariant K\"ahler cone structures on $\mathcal A$, is given by 
the space $\mathcal R$ of $T$-invariant smooth positive functions $r : \mathcal A \to
\mathbf R$ such that $i\partial\barpartial r^2$ is positive $(1,1)$-form:
$$ \mathcal R := \{ r : \mathcal A \to \mathbf R\ |\ T\text{-invariant},\ i\partial\barpartial r^2 > 0\}.$$
Since the Reeb vector field $\tilde \xi = Jr\frac{\p}{\partial r}$ is the real part of a holomorphic Killing vector field and $T$ is the maximal torus in $G$, 
$Jr\frac{\p}{\partial r}$ is in $\mathrm{Lie}(T)$ for any $r \in \mathcal R$.
The set of all Reeb vector fields corresponding to $r \in \mathcal R$ is the dual cone of the cone obtained as the moment map image of
$C(S)$, and is called the {\it Sasaki cone}.
We define the volume functional $\mathrm{Vol} : \mathcal R \to \mathbf R$ by
\begin{equation}\label{vol1}
 \mathrm{Vol}(r) = \mathrm{vol}(S_r)
 \end{equation}
where $\mathrm{vol}(S_r)$ denotes the volume of the Sasakian manifold $S_r = \{r = 1\}$ for $r \in \mathcal R$.
Let $r(t)_{-\epsilon < t < \epsilon}$ be a one parameter family in $\mathcal R$ with $r(0) = r$, and put
$Y := \frac{d}{dt}\vert_{t=0} \tilde \xi(t)$ where $\tilde \xi(t) = Jr(t)\frac{\p}{\partial r(t)}$.
Then the first variaton of $ \mathrm{Vol}(r)$ is given by
\begin{equation}\label{vol2}
\frac{d \mathrm{Vol}(r(t)) }{dt}\vert_{t=0} = -4(m+1)\int_{S_r} \eta(Y) dvol_r
 \end{equation}
where $dvol_r$ is the volume element of $S_r$, see \cite{FOW}, Proposition 8.3, or \cite{MSY2}, Appendix C1.
The second variation is given by 
\begin{equation}\label{vol3}
\frac{d}{dt}\vert_{t=0}\left( -4(m+1)\int_{S_{r(t)}} \eta(X) dvol_{r(t)}\right) = 4(m+1)(2m+4)\int_{S_r} \eta(X)\eta(Y) dvol_r,
 \end{equation}
 see \cite{FOW}, Proposition 8.4, or \cite{MSY2}, Appendix C2. 
 The second variation formula shows that the volume functional is convex.
 
A Sasakian manifold $S$ is called a Sasaki-Einstein manifold if it is an Einstein manifold as a Riemannian manifold. This occurs exactly when
$C(S)$ is a Ricci-flat K\"ahler cone (i.e. Calabi-Yau cone). From the view point of the K\"ahler sandwich, this occurs exactly
when the transverse K\"ahler structure of the Reeb flow is K\"ahler-Einstein with positive scalar curvature. 
A typical example is the $(2m+1)$-dimensional standard sphere which is Sasaki-Einstein. In this case, the cone is $\mathbf C^{m+1}$
which is Ricci-flat K\"ahler, the Reeb flow is the standard $S^1$-action, and the orbit space is the complex projective space which 
is a K\"ahler-Einstein manifold of positive scalar curvature.

When the Reeb flow generates an $S^1$-action the quotient space is a Fano orbifold. For general Sasakian structures the complex geometry of the local orbit spaces are described as ``basic'' geometry. For example, we have the basic $\partial$-operator $\partial_B$, 
the basic 
$\barpartial$-operator $\barpartial_B$, the basic 
Dolbeault cohomology $H^\ast_{\barpartial_B}$, the basic K\"ahler metric $g_B$, the basic K\"ahler form $\omega_B$, and the basic Ricci form $\rho_B$, the basic first Chern class $c_1^B$ and etc.
With these notations, the Sasaki-Einstein equation becomes 
$$ \rho_B = (2m+2) \omega_B.$$
Thus a necessary
condition for the existence of a Sasaki-Einstein metric is that the basic first Chern class is represented by a positive multiple 
of the basic K\"ahler class:
$$ 2\pi c_1^B = (2m+2) [\omega_B]$$
in $H^2_{\barpartial_B}$.
This last condition is equivalent to the topological condition that $c_1(D) = 0$ and that $c_1^B > 0$
where $D$ denotes the contact structure determined by the contact form $\eta$, see \cite{FOW}, Proposition 4.3. We say in this paper
that $S$ is transversely Fano if $c_1(D) = 0$ and $c_1^B > 0$. 
 
Let $\xi$ be the Reeb vector field on a Sasakian manifold $S$. A smooth function $f$ on $S$ 
is said to be basic if $\xi f = 0$. A basic function is obtained locally by pulling back a smooth function on the local orbit space of the Reeb
flow. A holomorphic vector field $Y$ in $\mathrm{Lie}(G)$ descends to a complex vector field on $S$ and also to a complex vector field
on each local orbit space of the Reeb flow, both of which we also denote by the same
letter $Y$. Then $Y$ is written on the local orbit space of the Reeb flow, which is K\"ahler, as
\begin{equation}\label{grad} Y = \mathrm{grad}_{g_B}^\prime u = g_B^{i\barj}\frac{\partial u}{\partial \overline{z^j}}\frac{\partial}{\partial z^i} \end{equation}
where $z^1, \cdots, z^m$ are local holomorphic coordinates and $g_B$ is the transverse K\"ahler metric 
on the local orbit space of the Reeb flow. There is a real valued basic function $F_B$ such that 
\begin{equation}\label{Ricci}
\rho_B - (2m+2)\omega_B = i\partial_B\barpartial_B F_B.
\end{equation}
Just as in the case of Fano manifolds (c.f. \cite{futaki88}, Theorem 2.4.3), there is an isomorphism between $\mathrm{Lie}(G)$ and
the space $\Lambda_{2m+2}$ of eigenfunctions $u$ of the elliptic operator $\Delta^F_B $ defined by
\begin{equation}\label{eigen}
\Delta^F_B u := \Delta_B u - \nabla^i u \nabla_iF_B
\end{equation}
where $\Delta_B = \barpartial_B^\ast \barpartial_B$ is the transverse $\barpartial_B$-Laplacian and $\nabla$ denotes the Levi-Civita connection of the transverse K\"ahler structure, see \cite{FOW}, Theorem 5.1.
Noting $\eta(Y)$ in (\ref{vol2}) is basic, if $\eta(Y) = u$ in $\Lambda_{2m+2}$, then
the right hand side of (\ref{vol2}) is equal to
\begin{eqnarray}
-2 \int_S (2m+2)u\ dvol &=& -2 \int_S (\Delta_B u - \nabla^i u \nabla_iF_B) dvol\\
&=& 2\int_S (\grad_{g_B}^\prime u)F_B\ dvl.
\end{eqnarray}
The right hand side is equal to $\mathrm{Fut}_\xi$  where $\xi$ is the Reeb vector field which is determined by the Sasakian structure
of $S$. This proves the volume minimization formula (\ref{derivative2}).

A Sasakian manifold 
$(S, g)$ is said to be toric if the K\"ahler cone manifold $C(S)$ is toric, namely 
$\dim_{\mathbf C} G = m+1$. When $S$ is toric and transversely Fano, Martelli-Sparks-Yau \cite{MSY2} 
showed that the volume functional is proper on the space $\Sigma$ of Reeb vector fields
of charge $n$, which is a slice in the the Sasaki cone, i.e. the dual cone
of the moment map image of $C(S)$. Since the volume functional is convex by (\ref{vol3}), there is a unique
minimum on $\Sigma$ at which $\mathrm{Fut}_\xi$ vanishes. In \cite{FOW} it is shown that for this minimum
$\xi$ there is a Sasaki-Einstein metric.
Uniqueness assertion is also shown in \cite{CFO}. To sum up the following holds.
\begin{thm}[\cite{FOW}, \cite{CFO}]\label{Main1}  Let $(S, g)$ be a compact toric Sasakian manifold with $c_1^B > 0$
and $c_1(D) = 0$. Then there exists a Sasaki-Einstein metric. Further, 
the identity component of the automorphism
group for the transverse holomorphic structure acts transitively 
on the space of all Sasaki-Einstein metrics. 
\end{thm}

In K\"ahler geometry, the Yau-Tian-Donaldson conjecture relates the existence problem of constant scalar curvature
K\"ahler (cscK for short)  metrics to K-stablity. Simlilarly in Sasakian geometry, the existence problem of constant scalar curvature
Sasaki (cscS for short)  metrics is related to K-stablity, see \cite{CollinsSzeke12}, \cite{CollinsSzeke15}, \cite{TiplervanCoev15}, \cite{BoyervanCoev2016} for example.

The cscS metrics are critical points of the Einstein-Hilbert functional $H : \mathcal R \to \mathbf R$ defined by
\begin{equation}\label{EH}
H(r) = \frac{\mathrm{TS}(r)^{m+1}}{\mathrm{Vol}(r)^m}
\end{equation}
where $\mathrm{TS}(r)$ denotes the total scalar curvature of $S_r$.
In the transversely Fano case, $\mathrm{TS}(r) = \mathrm{Vol}(r)$ and the Einstein-Hilbert functional
coincides with the volume functional. 
For general Sasakian manifolds, i.e. for Sasakian manifolds which are not necessarily transversely Fano,
it is known that the convexity fails for the Einstein-Hilbert functional, and there can be several critical
points, see Legendre \cite{Legendre11_2}, and also \cite{BHLT17}.
This fact has resemblance in the study of Einstein-Maxwell K\"ahler metrics  as can be seen in 
the ambitoric examples by LeBrun \cite{L2} on the one-point-blow-up of $\mathbf C\mathbf P^2$.
But it is shown by Boyer-Huang-Legendre \cite{BHL17} that all of the volume functional, the total scalar
curvature and the Einstein-Hilbert functional are proper in that they tend to $+\infty$ as the Reeb vector
field tends to the boundary of the Sasaki cone. This was shown by using the Duistermaat-Heckman formula.

The idea of volume minimization for Sasaki-Einstein metrics has been extended and applied to
algebraic geomrty. Odaka \cite{Odaka15} considered generalizations of the normalized volume
functional 
and Donaldson-Futaki invariant obtained as the derivative of the volume functional.
Odaka observed the decrease of the Donaldson-Futaki invariant along the minimal model 
program using the concavity of the volume functional.
 Li \cite{LiChi15_1}, \cite{LiChi15_2} considered normalized volume functional on the space of valuations
 on Fano manifolds
 and characterized K-semistabilty in terms of volume minimization. Note that when
 a Sasakian manifold is the circle bundle of an ample line bundle $L$ over $M$,
 then the Reeb vector field defines a valuation of the ring $\oplus_{k=0}^\infty H^0(M,L^k)$.
 In view of this, to define volume functional for valuations is natural. The normalization corresponds
 to the restriction of the Reeb vector fields to the ones with charge $n$. 
 On the other hand the Gromov-Hausdorff limit of a sequence of K\"ahler-Einstein manifolds is
 homeomorphic to a normal algebraic variety and admits a K\"ahler-Einstein metric in the sense 
 of pluripotential theory \cite{DonaldsonSun14}. The tangent cone at a singular point admits a Ricci-flat cone structure,
 and thus it is a cone over a Sasaki-Einstein manifold on the regular set. Li-Xu \cite{LiXu17}
 applies the volume minimization to show an algebraic nature of those tangent cones, answering
 to a question of Donaldson-Sun \cite{DonaldsonSun15}.
 See also \cite{LiLiu16}, \cite{LiXu16}.

\end{document}